\magnification=1200  
\overfullrule=0pt
\baselineskip=14pt

\font\tenbb=msbm10
\font\sevenbb=msbm7
\font\fivebb=msbm5
\newfam\bbfam
\textfont\bbfam=\tenbb \scriptfont\bbfam=\sevenbb
\scriptscriptfont\bbfam=\fivebb
\def\bb{\fam\bbfam}

\def\Cb{{\bb C}}
\def\Nb{{\bb N}}
\def\Pb{{\bb P}}
\def\Qb{{\bb Q}}
\def\Rb{{\bb R}}
\def\Zb{{\bb Z}}

\def\Oc{{\cal O}}

\def\a{\alpha}
\def\D{\Delta}
\def\lb{\lambda}
\def\om{\omega}
\def\s{\sigma}
\def\Si{\Sigma}
\def\t{\theta}
\def\ve{\varepsilon}

\def\ot{\otimes}
\def\ov{\overline}
\def\sbs{\subset}
\def\sm{\simeq}
\def\ts{\times}
\def\wh{\widehat}
\def\wt{\widetilde}

\def\ra{\rightarrow}
\def\hra{\hookrightarrow}
\def\longra{\longrightarrow}

\def\build#1_#2^#3{\mathrel{
\mathop{\kern 0pt#1}\limits_{#2}^{#3}}}

\def\hfl#1#2{\smash{\mathop{\hbox to 6mm{\rightarrowfill}}
\limits^{\scriptstyle#1}_{\scriptstyle#2}}}

\def\ldisplaylinesno #1{\displ@y\halign{
\hbox to\displaywidth{$\@lign\hfil\displaystyle##\hfil$}&
\kern-\displaywidth\rlap{$##$}
\tabskip\displaywidth\crcr#1\crcr}}
\catcode`\@=12


\centerline{\bf Secant varieties and successive minima}

\medskip

\centerline{\bf Christophe Soul\'e}

\smallskip

\centerline{\it CNRS et IHES}

\medskip

\centerline{10.10.2001}

\vglue 1cm

Let $X$ be a semi-stable arithmetic surface over the spectrum $S$ 
of the ring of integers in a number field $K$. We assume that the 
generic fiber $X_K$ is geometrically irreducible and has positive 
genus. Consider a line bundle $L$ on $X$, which is non-negative on 
every vertical fiber and equipped with an admissible metric 
(in the sense of Arakelov $[A]$) of 
positive curvature. In [S] we have shown that any non-trivial 
extension class of $L$ by the trivial line bundle is such that its 
$L^2$-norm is bounded away from zero.

\smallskip

In this paper we refine this analysis by showing that all 
successive minima of the euclidean lattice ${\rm Ext} \, (L, 
\Oc_X)$ admit an explicit lower bound, which can be computed in 
terms of the Arakelov intersection theory on $X$ (Th.~4). By Serre 
duality, this result gives also an upper bound for the successive 
minima of the global sections over $X$ of the tensor product of 
$L$ with the relative dualizing sheaf $\om_{X/S}$ of $X$ over $S$.

\smallskip

Our method to prove this result is a continuation of [S], Th.~2. 
The new ingredient is a study of the {\it secant varieties} 
$\Si_d$, $d > 0$, of the curve $X_K$,
when this curve is embedded into a projective 
space by the sections of $L \ot \om_{X/S}$. We consider the maximal 
dimension of a projective subspace contained in $\Si_d$. When 
$\deg \, (L) > 2d+2$, a result due to C.~Voisin asserts that this 
dimension is as small as possible (equal to $d-1$, Th.~1). 

\smallskip

The study of the projective subspaces of the secant variety 
$\Si_d$ is related to our initial question for the following 
reasons. On one hand, the first $k$ successive minimal vectors of 
${\rm Ext} \, (L , \Oc_X)$ span a rank $k$ sublattice; if
$k$ is big enough, this sublattice is not contained
in the secant variety.
On the 
other hand, if we know that an extension class $e$ does not lie in 
$\Si_d (K)$ for a large value of $d$, the proof of [S], Th.~2, 
gives a large lower bound for the $L^2$-norm of $e$.

\smallskip

The paper is organized as follows. 
In Section~1 we compute the degree of 
the secant varieties and we give Voisin's result 
on secant varieties of curves (Th.~1). We deduce from it
 a geometric result on extension classes
in Theorem 2. In Section~2 we propose another approach to 
linear subspaces in secant varieties of curves.
Though weaker than Theorem 1, it might be of independent
interest, especially because of Theorem 3, due to
A.~Granville. Granville's result computes the maximal length
of a  sequence
of binomial coefficients, in a given line of Pascal's triangle,
which have a common divisor. The answer
 turns out to be related to the
distribution of primes among positive integers.

 We start dealing with 
arithmetic surfaces in Section~3. We prove a ``transfer theorem'' 
for successive minima of dual lattices (Prop.~4), and we extend a 
result of Szpiro on the minimal height of algebraic points on $X$ 
(Prop.~5). We then prove our main theorem on successive minima 
over $X$ (Th.~4). We conclude with some explicit estimates when 
$L$ is a high  of $\om_{X/S}$ (Cor. 1 and 2).

\medskip

I thank  both C.~Voisin and A.~Granville
for discovering Th.~1 and Th.~3, respectively, when being
 asked for help. I also benefited from the advices 
of J.~Harris and O.~Ramar\'e.

\vglue 1cm

\noindent {\bf 1. Secant varieties}

\smallskip

\noindent {\bf 1.1.} Let $K$ be an algebraically closed field of 
characteristic zero, $C$ a smooth projective connected curve of 
genus $g$ over $K$, $\om_C$ the canonical line bundle on $C$, and 
$L$ a line bundle of degree $m = \deg \, (L) > 2$ over $C$. Denote 
by
$$
\Pb = \Pb (H^0 (C , L \ot \om_C)^*)
$$
the projective space defined by sections of $L \ot \om_C$ on $C$, and let
$$
j : C \ra \Pb
$$
be the projective embedding defined by these sections. For any 
positive integer $d > 0$, we let $\Si_d \sbs \Pb$ be the secant 
varieties of $d$-uples of points on $\Pb$, i.e.
$$
\Si_d = \bigcup_{\deg (D) = d} \ov D \, ,
$$
where $\ov D \sbs \Pb$ is the linear span of the effective divisor 
$D$ of degree $d$ over $C$ ([ACGH], [B]). We first compute the 
degree of this projective variety.

\medskip

\noindent {\bf Proposition 1.} {\it Assume $2d \leq \dim \Pb + 1$. 
Then $\Si_d$ has dimension $2d-1$ and the degree of $\Si_d$ in 
$\Pb$ is}
$$
\deg \, (\Si_d) = D (g , m , d) = \sum_{\a = 0}^{{\rm Inf} (d,g)} 
\left( \matrix{m+g-1-d-\a \cr d-\a \cr} \right) \left( \matrix{g 
\cr \a \cr} \right) \, .
$$

\noindent {\bf 1.2. Proof of Proposition $1$.} Let $M = L \ot 
\om_C$, let $C_d$ be the symmetric product of $d$ copies of $C$ and 
let $\D \sbs C \ts C_d$ be the universal divisor ([ACGH] p.~339). 
If $\pi_1 : C \ts C_d \ra C$ and $\pi_2 : C \ts C_d \ra C_d$ are 
the projections, we let
$$
E_L = \pi_{2*} (\Oc_{\D} \ot \pi_1^* \, M) \, .
$$
This is a vector bundle of rank $d$ on $C_d$ and the canonical map
$$
H^0 (C,M) \ot \Oc_{C_d} \ra E_L
$$
induces a morphism
$$
\a : \Pb \, (E_L^*) \ra \Pb
$$
the image of which is $\Si_d$ (op.cit. p.~340).

\smallskip

We first notice that $\a$ is a  birational isomorphism onto 
$\Si_d$. Indeed, by the general position theorem, [ACGH] p.~109, a 
general hyperplane $H$ in $\Pb$ cuts $C$ in exactly $\deg \, (M) = 
m + 2g - 2$ points, any $N$ of which are independent, where
$$
N := \dim \, (\Pb) + 1 = m+g-1 \, .
$$
In particular, $d$ points in $H \cap C$ span a linear subspace $Q 
\sbs \Pb$ of dimension $d-1$, which meets $C$ in exactly $d$ 
points. A general point on $Q$ has thus only one preimage by $\a$.

\smallskip

It follows that $\Si_d$ has dimension $2d-1$ and the degree of 
$\Si_d$ in $\Pb$ coincides with the integer $\deg \, (\xi^{2d-1})$ 
where
$$
\xi = c_1 (\a^* (\Oc (1))) \in {\rm CH}^1 (\Pb \, (E_L^*)) \, .
$$
The direct image of $\xi^{2d-1}$ by the projection
$$
\Pb \, (E_L^*) \ra C_d
$$
is the top Segre class of $E_L^*$ on $C_d$. Therefore $\deg 
\, (\xi^{2d-1})$ is the coefficient of $t^d$ in the formal power 
series in one variable
$$
s_t (E_L^*) = c_t \, (E_L^*)^{-1} \, .
$$
The total Chern series
$$
c_t (E_L^*) = \sum_{i \geq 0} \, t^i \, c_i (E_L^*) = c_{-t} (E_L)
$$
is computed in [ACGH], Lemma~2.5, p.~340. Fix a point $P_0 \in 
C(K)$ and denote by $x$ the class of $P_0 + C_{d-1}$ in ${\rm CH}^1 
(C_d)$, and by $\t \in {\rm CH}^1 (C_d)$ the restriction to $C_d$ 
of the class of the theta divisor on the Jacobian of $C$. We know 
from loc.cit. that
$$
c_t (E_L^*) = (1 + xt)^{-A} \, e^{-t\t / (1+xt)}
$$
where
$$
A = -d + \deg \, (L \ot \om_X) -g + 1 = m + g - 1 - d > 0 \, .
$$
Therefore
$$
\deg \, (\Si_d) = [c_t \, (E_L^*)^{-1}]_{t^d} = \sum_{\a = 0}^{{\rm 
Inf} (d,g)} \left( \matrix{ A - \a \cr d - \a \cr} \right) \, 
x^{d-\a} \ {\t^{\a} \over \a !} = \sum_{\a = 0}^{{\rm Inf} (d,g)} 
\left( \matrix{ A - \a \cr d - \a \cr} \right) \left( \matrix{ g 
\cr \a \cr} \right)
$$
by [ACGH], p.~343, last two lines. This proves Proposition~1. 
\hfill q.e.d.

\bigskip

\noindent {\bf 1.3.} By definition, $\Si_d$ is a union of linear 
subspaces of $\Pb$ of dimension $d-1$. We are interested in the 
dimensions of linear subspaces of $\Pb$ contained in $\Si_d$. 
\medskip

\noindent {\bf Theorem 1.} (C.~Voisin [V]) {\it Assume $m > 2d+2$. Then, there
is no linear subspace contained in  $\Si_d$ of dimension
bigger than $d-1$. Furthermore, the only linear subspaces
of dimension $d-1$ are those spanned by effective divisors of degree
$d$ on $C$.}

\medskip

 C.~Voisin conjectures that the hypothesis
 $m > 2d$ is enough for Theorem 1 to hold. 
 When $m = 2d+1$ or $m=2d+2$, we have a weaker bound 
than Theorem 1:

\smallskip

\medskip

\noindent {\bf Proposition 2.}

 {\it If $m = 2d+1$ or $m=2d+2$, any linear 
 subspace contained in  $\Si_d$ has dimension
 at most $2d-g-1$.}

\bigskip

\noindent {\bf 1.4. Proof of Proposition 2.} 
By induction on $d$, we may consider a  linear 
space $Q \sbs \Si_d$ such that $Q \not\sbs \Si_{d-1}$.
 The morphism
$$
\a : \Pb \, (E_L^*) \ra \Pb
$$
is an embedding when restricted to
 $\Si_d - \Si_{d-1}$ ([B] Lemma~1.2~a)). 
Therefore the closure $Y$ in $\Pb \, (E_L^*)$ of $\a^{-1} (Q - 
\Si_{d-1})$ is birationally equivalent to $Q$. 
In particular, its Albanese variety is trivial.
Let ${\rm Pic}^{(d)} 
(C)$ be the component of the Picard variety of $C$ parametrizing 
line bundles of degree $d$, and $\s : C_d \ra {\rm Pic}^{(d)} (C)$ 
the canonical map. The composite map
$$
Y \longra \Pb \, (E_L^*) \build \longra_{}^{\pi} C_d \build 
\longra_{}^{\s} {\rm Pic}^{(d)} (C)
$$
must be constant since $Alb(Y)$ is trivial. 
Let $x_0 \in {\rm Pic}^{(d)} (C)$ be the 
image of $Y$. We denote by $S$ (resp. $T$) 
 the inverse of $x_0$ by the map 
$\s$ (resp. $\s \circ \pi$). Since $\deg_K (M) > 2g-2$, $\s$ is a 
projective bundle, hence $S \sm \Pb^{d-g}$ ([ACGH], VII, 
Prop.~2.1), and the map $\pi : T \ra S$ is also a projective bundle,
of relative dimension $d-1$. It 
follows that the dimension of $Q$ is at most $2d-g-1$.

\hfill q.e.d.

\bigskip

\noindent {\bf 1.5.} We keep the notation of \S~1.1. Let $L^{-1}$ 
be the dual of $L$ and
$$
e \in {\rm Ext}^1 (L , \Oc_C) = H^1 (C , L^{-1}) 
$$
any extension class of $L$ by the trivial bundle. There exists a 
rank two vector bundle $E$ over $C$ and an exact sequence of 
coherent sheaves
$$
0 \ra \Oc_C \ra E \ra L \ra 0 \leqno (1)
$$
which is classified by $e$. This extension is uniquely determined 
by $e$ up to isomorphism.

\smallskip

Recall that $E$ is called semi-stable (resp. stable) when any 
invertible subsheaf $M \sbs E$ satisfies
$$
\deg \, (M) \leq \deg \, (E) / 2 = m/2
$$
(resp. $\deg \, (M) < m/2$).

\medskip

\noindent {\bf Theorem 2.} {\it Let $d$ be a positive integer  
and $V$  a $K$-vector space of $ H^1 (C , 
L^{-1})$.

\noindent{\rm i)} If $m > 2d +2$ and
$
\dim_K (V) > d $
 there exists a non-trivial extension class 
$e \in V$ with the following 
property: given any line bundle
$M \sbs V$, one has
$$
\deg \, (M) < m-d \, .
$$

\noindent{\rm ii)}
When $m = 2d+1$ (resp. $m = 2d+2$) and 
$
\dim_K (V) > 2d-g $
there exists a non-trivial $e\in V$ such that
 $E$ is stable (resp. semi-stable).}

\medskip

\noindent {\bf 1.6. Proof of Theorem 2.} Let $E$ be a non-trivial 
extension of $L$ by ${\Oc}_C$ as in (1) and let $M \sbs E$ be a 
line bundle of positive degree. Then the map $M \ra L$ is 
non-trivial (otherwise $M \sbs \Oc_C$), therefore there exists an 
effective divisor $D$ on $C$ such that the map $M \ra L$ induces an 
isomorphism
$$
M \build \longra_{}^{\sim} L (-D) \, .
$$
When restricted to $L (-D)$ the extension (1) becomes trivial, 
i.e. $e$ lies in the kernel of the map
$$
H^1 (X,L^{-1}) \ra H^1 (X,L (-D)^{-1}) \, .
$$
The exact sequence of sheaves on $C$
$$
0 \ra L^{-1} \ra M^{-1} \ra M^{-1} \ot \Oc (D) \ra 0
$$
gives an exact sequence of cohomology groups
$$
H^0 (C , M^{-1} \ot \Oc (D)) \build \longra_{}^{\partial_D} H^1 
(C,L^{-1}) \longra H^1 (C,M^{-1}) \, .
$$
Therefore $e$ lies in the image of $\partial_D$. But the map 
$\partial_D$ coincides with the restriction of
$$
\a : \Pb \, (E_L^*) \ra \Pb
$$
to the fiber at $D \in C_d$ of the projection 
map $\Pb \, (E_L^*) \ra 
C_d$ (see [ACGH], p.~340, or [B], Observation~2, p.~451). In other 
words, if $d = \deg \, (D) = \deg \, (L) - \deg \, (M)$, the class 
of $e$ in $\Pb$ must lie in the secant variety $\Si_d$. By  
Theorem 1 we know that if $m>2d+2$ and $\dim (V) > d $, 
there exists a non-zero
$e \in V$ whose class in $\Pb = \Pb \, (H^1 (C , 
L^{-1}))$ (Serre duality) is not in $\Si_d$. This proves the 
first assertion in Theorem~2.

\smallskip

When $m = 2d+1$ (resp. $m = 2d+2$) we get
 from Proposition 1, under the hypothesis in ii), that
$
\deg \, (M) < d+1$, hence
$$ \ \deg \, (M) < \deg \, 
(E) / 2 
$$
(resp. $\deg \, (M) \leq d+1 = \deg \, (E) / 2$). Therefore $E$ is 
stable (resp. semi-stable).

\hfill q.e.d.

\vglue 1cm

\noindent {\bf 2. On the divisibility of binomial coefficients}

\smallskip

\noindent {\bf 2.1.} For any integer $n > 0$, let $b (n) \geq 0$ be 
the smallest integer $b$ such that the set of binomial coefficients 
$\left( \matrix{n \cr m \cr} \right)$, where $b < m < n-b$, has a 
(non-trivial) common divisor.

\medskip

\noindent {\bf Theorem 3.} (A. Granville) {\it The integer $b(n)$ is the 
smallest integer of the form $n - p^k$, where $p^k$ is a prime 
power less or equal to $n$.}

\medskip

\noindent {\bf 2.2.} Let $c(n)$ be the smallest integer of the form 
$n - p^k$, with $p^k \leq n$ a prime power. The inequality $b(n) 
\leq c(n)$ is not hard to prove. First, when $n = p^k$ is itself a 
prime power, it is clear that $b(n) = 0$. Indeed (we follow here a 
suggestion of Tamvakis), for any $m > 0$, with $m < p^k$, the 
number
$$
\left( \matrix{ p^k - 1 \cr m \cr} \right) = {(p^k - 1) (p^k - 2) 
\ldots (p^k - m) \over m!}
$$
is congruent to $(-1)^m$ modulo $p$. Therefore
$$
\left( \matrix{ p^k \cr m \cr} \right) = \left( \matrix{ p^k - 1 
\cr m \cr} \right) + \left( \matrix{ p^k - 1 \cr m - 1 \cr} \right)
$$
is divisible by $p$.

\smallskip

Now, given any integer $n > 1$, we can find a prime power $p^k$ 
such that
$$
{n \over 2} < p^k \leq n
$$
(this fact is due to 
Chebyshev, see (6) below for a stronger statement).
Let $m$ be any 
integer such that
$$
n - p^k < m \leq p^k \, . 
$$
We have
$$
\left( \matrix{ n \cr m \cr} \right) = \left( \matrix{ n - 1 \cr m 
\cr} \right) + \left( \matrix{ n - 1 \cr m - 1 \cr} \right) \, . \leqno 
(2)
$$
If $n$ is not a prime power, both $m$ and $m-1$ are in the interval 
$] (n-1) - p^k , p^k ]$i and $p^k \leq n - 1$.
 If we fix the prime power $p^k$ and 
proceed by induction on $n$, with $p^k \leq n < 2 p^k$ (starting 
with the case $n = p^k$ treated above)
we conclude that both $\left( 
\matrix{ n-1 \cr m \cr} \right)$ and $\left( \matrix{ n-1 \cr m-1 
\cr} \right)$ are divisible by $p$, therefore $p$ divides $\left( 
\matrix{ n \cr m \cr} \right)$ as soon as $n - p^k < m \leq p^k$. 
This implies
$$
b (n) \leq n - p^k \, .
$$
Therefore $b(n) \leq c(n)$.

\medskip

\noindent {\bf 2.3.} It is more difficult to check that $b(n) = 
c(n)$. For each prime $p$ let $b_p$ be the largest integer $b \leq 
n/2$ such that $p$ does not divide $\left( \matrix{ n \cr b \cr} 
\right)$. We have
$$
b (n) = \ \build {\rm Inf}_{p \, {\rm a \, prime}}^{} \ b_p \, .
$$
Now write $n$ and $m$ in base $p$:
$$
n = n_k \, p^k + n_{k-1} \, p^{k-1} + \cdots + n_0
$$
$$
m = m_k \, p^k + m_{k-1} \, p^{k-1} + \cdots + m_0
$$
with $0 \leq m_i$, $n_i \leq p-1$ and $n_k > 0$ (we assume $m \leq 
n$). Kummer showed that $p$ divides $\left( \matrix{ n \cr m \cr} 
\right)$ if and only if $n_i < m_i$ for some $i$
([R  ] p.23-24). In 
particular, given any integer $r$ such that $0 \leq r \leq n_k$, 
$p$ does not divide $\left( \matrix{ n \cr r \, p^k \cr} \right)$.

\smallskip

Assume $n_k \geq 2$ and let $r$ be the integral part of $n_k / 2$. 
We have
$$
r \geq 1 \quad \hbox{and} \quad r \geq (n_k - 1) / 2 \, ,
$$
therefore
$$
n < (n_k + 1) \, p^k \leq 2 (r+1) \, p^k \leq 4 r \, p^k \, .
$$
On the other hand
$$
b_p \geq r \, p^k \, .
$$
Therefore
$$
b_p > n/4 \, . \leqno (3)
$$

As we shall see below (cf. (6)),
 for any integer $n \geq  2$ there exists a prime 
$p$ such that $3n / 4 \leq p \leq n$. By Kummer's criterion 
recalled above, we know that $p$ divides $\left( \matrix{ n \cr m 
\cr} \right)$ when $n-p < m < p$. Therefore
$$
b_p \leq n - p \leq n/4 \, ,
$$
which contradicts (3).

\smallskip

So, when computing $b_p$, we can assume that $n_k = 1$. This 
implies that $n < 2 \, p^k$. When $n - p^k < m < p^k$ there exists 
$i < k$ with $m_i > n_i$. Therefore, by Kummer's criterion (or by 
2.2), $p$ divides $\left( \matrix{ n \cr m \cr} \right)$. On the 
other hand $p$ does not divide $\left( \matrix{ n \cr n - p^k \cr} 
\right)$, again by Kummer's criterion. Therefore
$$
b_p = n - p^k \leqno (4)
$$
hence $b(n) = c(n)$. This ends the proof of Theorem~1. \hfill 
q.e.d.

\bigskip

\noindent {\bf 2.4.} Given two real valued 
functions $f(n)$ 
and $g(n)$
on positive integers (or real numbers),
we shall write $f(n) \ll g(n)$ when there exists a positive constant
$C>0$ such that, for all $n$, $f(n) \leq C g(n)$.

According to [ B-H ], there exists a positive 
constant $C$ such that, for any integer $n > 1$, there is a 
prime number $p$ with $p \leq n$ and 
$$
p \geq n - C \, n^{0.535} \, .
$$
This implies
$$
b(n) \ll \,  n^{0.535} \, . \leqno (5)
$$
The Riemann hypothesis is known to imply
$$
b (n) \ll \,  \sqrt n \ \log \, (n) 
$$
([ I ] Th. 12.10).

\smallskip

 It was shown by Nagura [N] that if $x\geq 25$
 there exists a prime $p$ with
 $$x < p < 6x/5 \, .$$
 Given any integer $n \geq 30$, write $n = 6y + r$
 with $0 \leq r < 6$. If we apply the previous inequality
 to $x = 5y$ we get that there exists a prime $p$
 such that
 $$ (5n/6)-4 \leq p < n \, . $$
 This implies in particular
$$
b(n) \leq n/4 \leqno (6)
$$
when $n \geq 30$. The inequality
(6) can  be checked directly when 
$2 \leq n \leq 30$.

\bigskip

\noindent {\bf 2.5.} From Theorem~3 one can also derive 
estimates on the sum of the first $n$ values of the function 
$b$ when $n$ gets large. Indeed, if $p_1 , p_2 , p_3 , 
\ldots$ is the list of prime numbers in increasing order, 
one has
$$
\matrix{
\displaystyle \sum_{j=1}^n b(j) &\leq &\displaystyle 
\sum_{p_k \leq n} \ \sum_{p_k+1}^{p_{k+1}} b(j) \hfill \cr 
\cr
&\leq &\displaystyle \sum_{p_k \leq n} \left( 
\sum_{j=0}^{p_{k+1} - p_k} j \right) \hfill \cr \cr
&= &\displaystyle \sum_{p_k \leq n} (p_{k+1} - p_k) (p_{k+1} 
- p_k - 1) / 2 \hfill \cr \cr
&\leq &\displaystyle {1 \over 2} \sum_{p_k \leq n} (p_{k+1} 
- p_k)^2 \hfill \cr \cr
&\ll_{\ve} &n^{23/18+\ve} \, , \hfill \cr
}
$$
where the symbol $\ll_{\ve}$ means that the 
constant involved depends on $\ve$ (but not on $n$),
and the last inequality is a result due to 
D.R.~Heath-Brown ([I] Th. 12.17, (12.117)). So we 
have
$$
\sum_{j=1}^n b(j) \ll_{\ve} n^{23/18+\ve} \, . \leqno (7)
$$
A result of Selberg (cf. [I] p.~349) says that the Riemann 
hypothesis implies that
$$
\sum_{j=1}^n b(j) \ll n \log (n)^3 \, .
$$

\bigskip

\noindent {\bf 2.6.}
The following result is weaker than C.Voisin's Theorem 1. Still, because
of the facts described in \S 2.1 to 2.5, it puts strong constraints
on the dimension of linear subspaces contained in secant varieties.
We use the notation of \S 1.1.

\medskip

\noindent {\bf Proposition 3}
{\it If $m > 2d$ and
$$
b (m+g-1-d) \leq m+2g-1-2d \, ,
$$
then
$$
\ve (d) < b (m+g-1-d) \, .
$$
}

\noindent {\bf 2.7. Proof of Proposition 3.} From Theorem~3 we get
$$
b(n) \leq b (n-1) + 1
$$
therefore the quantity $b (m+g-1-d) +d-1$ is an increasing function 
of $d$. Therefore, by induction on $d$, 
to prove Proposition 3, we may consider a linear 
space $Q \sbs \Si_d$ such that $Q \not\sbs \Si_{d-1}$ and $\dim Q = 
d-1+\ve (d)$. As in the proof of Proposition 1,
we let  $Y$
be the closure in $\Pb \, (E_L^*)$ of $\a^{-1} (Q - 
\Si_{d-1})$ and
 $x_0 \in {\rm Pic}^{(d)} (C)$ be the 
image of $Y$. Denote by $S$ (resp. $T$) 
 the inverse of $x_0$ by the map 
$\s$ (resp. $\s \circ \pi$). Since $\deg_K (M) > 2g-2$, $\s$ is a 
projective bundle, hence $S \sm \Pb^{d-g}$ ([ACGH], VII, 
Prop.~2.1). Let $x \in {\rm CH}^1 (S)$ be the class of an 
hyperplane. The map $\pi : T \ra S$ is also a projective bundle,
of relative dimension $d-1$. We let $\xi_0 \in {\rm CH}^1 
(\Pb)$ be the class of an hyperplane, hence $\xi = \a^* (\xi_0) \in 
{\rm CH}^1 (T)$ is a generator of the ring ${\rm CH}^* (T)$ over 
${\rm CH}^* (S)$. In particular, $x$ and $\xi$ span the ring 
${\rm CH}^* (T)$. From [ACGH], Lemma~2.5, p.~340, one can compute 
the Segre classes of $E_L^*$ restricted to $S$. If we let
$$
s_t \, (E_L^*) = \sum_{i \geq 0} \, s_i \, (E_L^*) \, t^i
$$
we get (see \S~1.2)
$$
s_t \, (E_L^*)_{\mid S} = (1+xt)^{m+g-1-d} \, .
$$
Therefore, for any $i \geq 0$,
$$
\pi_* (\xi^{d-1+i}) = s_i \, (E_L^*)_{\mid S}  = \left( \matrix{ 
m+g-1-d \cr i \cr} \right) \, x^i \, . \leqno (8)
$$
Since $\a$ restricted to $Y$ is a birational isomorphism of $Y$ 
with the linear subspace $Q$ in $\Pb$, we have
$$
[Y] \cdot \xi^{d-1+\ve (d)} = 1 \, , \leqno (9)
$$
where $[Y]$ is the class of $Y$ in the appropriate Chow group of 
$T$, namely ${\rm CH}^k (T)$ with
$$
k = 2d-1-g-(d-1+\ve (d)) = d-g-\ve (d) \, .
$$
Let us write $[Y]$ in terms of our basis of ${\rm CH}^* (T)$:
$$
[Y] = \sum_{i+j=d-g-\ve (d)} \, n_i \, x^i \, \xi^j \, .
$$
From (8) and (9) we get 
$$
1 = \sum_{i+j=d-g-\ve (d)} \, n_i \, x^i \, \xi^{j+d-1+\ve (d)} = 
\sum_{i+j=d-g-\ve (d)} \, n_i \left( \matrix{m+g-1-d \cr j+\ve 
(d) \cr} \right) \, .
$$
Therefore the binomial coefficients $\left( \matrix{m+g-1-d \cr b 
\cr} \right)$ are prime to each other when
$$
\ve (d) \leq b \leq d-g \, .
$$
Let $A = m+g-1-d$ and $n = A-d+g = m+2g-1-2d$. By the definition of 
$b(A)$ there are two possibilities:

\smallskip

\item{i)} $n < \ve (d)$ and $n < b(A)$

\smallskip

\item{ii)} $\ve (d) \leq n$ and $\ve (d) < b(A)$.

\smallskip

\noindent Therefore, if $b(A) \leq n$ we must have
$$
\ve (d) < b(A) \, .
$$
\hfill q.e.d.

\noindent {\bf 3. Successive minima on arithmetic surfaces}

\smallskip

\noindent {\bf 3.1.} Let $K$ be a number field, $\Oc_K$ its ring of 
integers and $S = {\rm Spec} \, (\Oc_K)$. Consider an hermitian 
vector bundle $\ov V$ of rank $N+1$ over $S$, i.e. a finitely 
generated projective $\Oc_K$-module $V$ and, for every complex 
embedding $\s : K \ra \Cb$, an hermitian scalar product on the 
complex vector space $V_{\s} = V \build \ot_{\Oc_K}^{} \Cb$ defined 
by $\s$; furthermore these scalar products are invariant by complex 
conjugation.
We denote by $\wh{\deg} (\ov V) \in \Rb$ the arithmetic degree of 
$\ov V$ (see [BGS] (2.1.11), (2.1.15)).

\smallskip

Given any positive integer $p$, $1 \leq p \leq N+1$, we let $\lb_p 
(\ov V)$ be the $p$-th minimum of $\ov V$, i.e. the infimum of the 
set of real numbers $\lb$ such that there exist $p$ vectors $e_i 
\in V$, $1 \leq i \leq p$, which are linearly independent in $V 
\build \ot_{\Oc_K}^{} K$ and such that, for any $i$ and any complex 
embedding
$\s : K \hra \Cb$, $\log \, \Vert e_i \Vert_{\s} \leq \lb$.
Furthermore, let $\ell_p (\ov V)$ be the minimal
projective height of $\Pb (F)\subset \Pb (V)$
when $F$ is a subbundle of rank $p$ in $V$. 

\smallskip

Let $V^* = {\rm Hom} \, (V , \Oc_K)$ be the dual of the 
module $V$. We equip $V^*$ with the metric dual to the one chosen 
on $V$.
Finally, if $B_n$ is the euclidean volume
of the unit ball in $\Rb^n$, let 
$r_1$ (resp. $r_2$) be the number of real (resp. complex) places of 
$K$, and let $\D_K$ be the absolute discriminant of $K$.
We introduce the following constant 
$$
C (N,K) = (N+1) (r_1 + r_2) \, \log (2) + (N+1) (\log \vert \D_K 
\vert ) / 2 - r_1 \, \log \, B_N - r_2 \, \log \, B_{2N+2}  
\, .
$$

\noindent {\bf Proposition 4.} {\it Given any $p$ with $1 \leq p 
\leq N+1$, the following inequalities hold:}
$$
  \ell_{N+1-p} (\ov{V}^*)\leq  [K:\Qb] \,\sum_{j=1}^{p} \, 
\lb_j (\ov V) \leq  C (N,K) + \ell_{N+1-p} (\ov{V}^*)
 \, .
$$

\noindent {\bf 3.2. Proof of Proposition 4.}
Since
$\wh{\deg} (\ov{V}^*) = - \wh{\deg} (\ov V)$, it follows from [BGS]
 (4.1.3) that
$$
\ell_p (\ov V) = \ell_{N+1-p} (\ov{V}^*) + \wh{\deg} (\ov V) \, .
$$
From [BGS] Theorem 5.2.4, this implies
$$\ell_{N+1-p} (\ov{V}^*) =  \ell_p (\ov V) - \wh{\deg} (\ov V)
\leq [K:\Qb] \, \sum_{j=1}^{p} \lb_j (\ov V) \, .$$
From [BGS] (5.2.14) and (5.2.15) (where $\lb_j (\ov V)$ is denoted
$\lb'_j$), and from the result of Bombieri and Vaaler [Bo-Va], we
get
$$
\matrix{
\displaystyle [K:\Qb] \, \sum_{j=1}^{p} \, \lb_j (\ov V) &\leq
&\displaystyle C (N,K) - \wh{\deg} \, (\ov V) - [K:\Qb] \,
\sum_{j=p+1}^{N+1} \, \lb_j (\ov V) \hfill \cr
\cr
 &\leq &\displaystyle C (N,K) - \wh{\deg} \, (\ov V) + \ell_p (\ov
 V) = C (N,K) + \ell_{N+1-p} (\ov{V}^*) \, . \cr
 } 
 $$
This concludes the proof of Proposition 4. \hfill q.e.d.

\bigskip

\noindent {\bf 3.3.} We keep the notation of section~3.1. Let $\Pb 
= {\rm Proj} \, ({\rm Sym} \, (V^*))$ be the projective space of $V$ 
and $\Pb_K$ its generic fiber. Consider a closed subvariety $\Si 
\sbs \Pb_K$ of projective degree $D = \deg \, (\Si)$. Let
$$
\lb_{\max} (\ov V) = \lb_{N+1} (\ov V)
$$
be the last successive minimum of $\ov V$.

\medskip

\noindent {\bf Proposition 5.} {\it There exists a non zero vector 
$v \in V$ such that, for all $\s : K \hra \Cb$, 
$$
\log \, \Vert v \Vert_{\s} \leq \lb_{\max} (\ov V) + \log \, (D 
(N+1)) \, ,
$$
and such that the point $[v] \in \Pb (K)$ does not lie in $\Si 
(K)$.}

\medskip

\noindent {\bf Proof of Proposition 5.} Choose a basis $\{ e_1 , 
\ldots , e_{N+1} \}$ of $V \build \ot_{{\Oc}_K}^{} K$ such that 
$e_i \in V$ and
$$
\log \, \Vert e_i \Vert = \lb_i
$$
for all $i=1 , \ldots , N+1$. Using this basis, we can identify 
$\Pb$ with $\Pb^N$.

\smallskip

We first notice that there exists an homogeneous polynomial $F (X_1 
, \ldots , X_{N+1})$ of degree $D$ such that $\Si$ lies in the zero 
set of $F$. Indeed we can find a linear projection $\pi : \Pb^N 
\cdots\!\!\ra \Pb^M$ such that $\pi (\Si)$ is an hypersurface of 
degree $D$ in $\Pb^M$. If $\wt F = 0$ is an equation for $\pi 
(\Si)$, we can take $F = \wt F \circ \pi$.

\smallskip

Now let $\Phi$ be the finite set of vectors of the form 
$$
v = \sum_i \, n_i \, e_i
$$
where $n_i \in \Zb$ and $0 \leq n_i \leq D$ for all $i = 1 , 
\ldots , N+1$. Assume that $F(v) = 0$ for any $v \in \Phi$. Let us 
write
$$
F(X) = \sum_{\a} \, r_{\a} \, X^{\a} \, ,
$$
where $\a$ runs over multi-indices $\a = (\a_1 , \ldots , \a_{N+1}) 
\in \Nb^{N+1}$ of degree $\a_1 + \cdots + \a_{N+1} = D$ and $X^{\a} 
= X_1^{\a_1} \ldots X_{N+1}^{\a_{N+1}}$. We get
$$
\sum_{\a} \, r_{\a} \, n_{\a} = 0 \, , \leqno (10)
$$
where $n_{\a} = n_1^{\a_1} \ldots n_{N+1}^{\a_{N+1}}$ for all
$(n_i) \in \Zb^{N+1}$ such that $0 \leq  n_i \leq D$ for all $i = 1, 
\ldots , N+1$. We can view (10) as a system of linear equations in the 
variables $r_{\a}$, with coefficients $n_{\a}$. Its determinant is the 
determinant of the tensor product of $N+1$ Vandermonde 
$(D+1)\times (D+1)$ square 
matrices of the form $(n^j)$, with $j = 0,1,\ldots , D$ and $0 
\leq n \leq D$. Up to sign, this determinant is
$$
\prod_{n \ne m} \, (n-m)^{N+1} \, ,
$$
where $n,m$ run over all integers such that $0 \leq n , m \leq D$. 
Since this determinant does not vanish we get a contradiction.

\smallskip

Therefore there exists $v \in \Phi$ such that $F(v) \ne 0$, hence 
$[v] \notin \Si$. Now, for all $\s : K \hra \Cb$, we have
$$
\log \, \Vert v \Vert_{\s} = \log \, \left \Vert \sum_i \, n_i \, 
e_i \right\Vert_{\s} \leq \max_i \log \, \Vert e_i \Vert + \log \, 
(N+1) + \log \, (D) = \lb_{\max} (\ov V) + \log \, (D (N+1) )  
\, .
$$
\hfill q.e.d.

\bigskip

\noindent {\bf 3.4.} Let $K$ be a number field and $S = {\rm Spec} 
\, (\Oc_K)$. From now on, $f : X \ra S$ will denote a semi-stable 
curve over $S$ with generic fiber a geometrically irreductible 
curve $X_K$ of genus $g \geq 2$. Let $\om_{X/S}$ be the relative 
dualizing sheaf, and $\mu$ the K\"ahler form introduced by Arakelov 
in [A], \S~4. We endow $\om_{X/S}$ with the same metric as in [A].

\smallskip

Let $\ov L$ be an hermitian line bundle on $X$ such that the 
curvature of $\ov L$ is a positive multiple of $\mu$, i.e.
$$
c_1 (\ov L) = m \, \mu
$$
where $m = \deg_K (L)$ is a positive integer. For any algebraic 
point $P$ on $X_K$, consider its normalized height $h_{\ov L} (P) 
\in \Rb$, which is defined as follows. Let $K'$ be the field of 
definition of $P$, $\Oc_{K'}$ its ring of integers, and $s : {\rm 
Spec} \, (\Oc_{K'}) \ra X$ the section defined by $P$ (using the 
valuative criterion for properness). Then
$$
h_{\ov L} (P) = \wh{\deg} \, (s^* \, \ov L) / [K' : K] \, .
$$
When $P$ varies, $h_{\ov L} (P)$ is bounded below and we let
$$
e (\ov L) = \ \build{\rm inf}_{P}^{} \, h_{\ov L} (P) \, .
$$

Given two hermitian line bundles ${\ov L}_1$ and ${\ov L}_2$ on $X$ 
we write
$$
{\ov L}_1 \cdot {\ov L}_2 = f_* (\wh{c}_1 ({\ov L}_1) \cdot \wh{c}_1 
({\ov L}_2)) \in \Rb
$$
the arithmetic intersection number of their first Chern classes. The 
following result is due to Szpiro when $\ov L = \ov{\om}_{X/S}$ 
([E], Th.5 b)).

\medskip

\noindent {\bf Proposition 6.} {\it The following inequality holds:}
$$
e (\ov L) \geq {g \, {\ov L}^2 \over 2 m} - {\ov L \cdot 
\ov{\om}_{X/S} \over 2} + {m \, \ov{\om}_{X/S}^2 \over 8g} \, . 
\leqno (11)
$$

\medskip

\noindent {\bf 3.5. Proof of Proposition 6.}  Let $P$ be an 
algebraic point on $X_K$, defined over $K'$. After extending scalars 
from $K$ to $K'$ (this multiplies both the height of
$P$ and the right hand side of (11) by $[K' : K]$),
we may assume that $K=K'$. Let $s : S \ra X$ be the 
section defined by $P$ and $D$ the Zariski closure of $\{ P \}$ in 
$X$. We endow $\Oc (D)$ with its canonical admissible metric $[A]$ 
and write $\ov D$ instead of $\ov{\Oc (D)}$. If we apply the Hodge 
index theorem as in [S], (33), we get the inequality
$$
{\ov L}^2 - 2m \, \ov L \, \ov D + m^2 \, {\ov D}^2 \leq 0 \, .
$$
The adjunction formula $[A]$
$$
{\ov D}^2 = - \ov{\om}_{X/S} \cdot \ov D
$$
implies
$$
{\ov L}^2 \leq (2m \, \ov L + m^2 \, \ov{\om}_{X/S}) \cdot \ov D
$$
i.e., if $\ov M = {\ov L}^{\ot 2} \ot \ov{\om}_{X/S}^{\ot m}$,
$$
\left( {\ov M - m \, \ov{\om}_{X/S} \over 2} \right)^2 \leq m \, \ov 
M \cdot \ov D \, .
$$
Since $m' = \deg_K (M) = 2mg$, we get
$$
h_{\ov M} (P) \geq {1 \over 2g \, m'} \left( {2g \, \ov M - m' \, 
\ov{\om}_{X/S} \over 2} \right)^2
$$
hence
$$
e (\ov M) \geq {g \, {\ov M}^2 \over 2m'} - {\ov M \cdot \ov{\om}_{X/S} 
\over 2} + {m' \, \ov{\om}_{X/S}^2 \over 8g} \, . \leqno (12)
$$
Now the inequality (11) is true for $\ov L$ if and only if it is 
true for a positive power of $\ov L$, and such a power can be written
$$
{\ov L}^{\ot k} = \ov{L'}^{\ot 2} \ot \ov{\om}_{X/S}^{\ot \deg_K 
(L')} \, ,
$$
where ${\ov L'}$ satisfies the same positivity assertion as $\ov L$.
 It follows from (12) that (11) is true.

\medskip

\noindent {\bf 3.6.} We keep the notation of \S~2.4 and we 
also assume that the degree of the restriction
of $L$ to any component of a fiber of $X$ over
$S$ is nonnegative. Consider 
the group
$$
H^1 (X , L^{-1}) = {\rm Ext} \, (L , \Oc_X)
$$
of extensions of $L$ by the trivial line bundle. This is a finitely 
generated projective $\Oc_{K}$-module. We equip $H^1 (X , L^{-1}) 
\ot_{\Zb} \, \Cb$ with the $L^2$-metric defined by $\mu$ and the chosen 
metric on $L$. The Serre duality isomorphism
$$
H^1 (X , L^{-1}) \sm H^0 (X , L \ot \om_{X/S})^*
$$
is compatible with the $L^2$-metrics on both cohomology groups. 
Since $H^0 (X , L^{-1}) = 0$, the Riemann-Roch theorem reads
$$
{\rm rk} \, H^1 (X, L^{-1}) = m + g - 1 \, .
$$
For any integer $k = 1 , \ldots , m+g-1$, we let $\lb_k (\ov L)$ 
 be the $k$-th successive minimum of the 
metrized $\Oc_K$-lattice $H^1 (X , L^{-1})$, and we let
$$\mu_k (\ov L ) =\ell_{m+g-1-k}(H^0 (X , L \ot 
\om_{X/S})) \, . $$

\medskip

\noindent {\bf Theorem 4.} {\it 

\noindent{\rm i)} Let $k > 1$ be such that $m > 2k$. 
Then the following inequalities hold:
$$
\lb_k (\ov L) + 1 \geq [k (\ov{L}^2 - 2 \, m \, e \, (\ov L)) + m^2 e 
\, (\ov L) - \log (D (g,m,k - 1) (m+g)) [K:\Qb]] / (m^2 [K:\Qb]) \, \leqno (13)
$$
and
$$
\mu_k (\ov L) \geq - C \, (m+g-2,K) \, + 
$$
$$
\left[ ( {k(k+1) \over 2} 
(\ov{L}^2 - 2 \, m \, e \, 
(\ov L)) + k \, m^2 e \, (\ov L) - \sum_{j=1}^k \log (D (g,m,j-1) 
(m+g)) [K:\Qb]
\right] / m^2  \, . \, \leqno (14)
$$
\noindent{\rm ii)} When $m$ is odd (resp. when $m$ is even
and $g>0$) we have
$$
\lb_{m-g-1} (\ov L) + 1 
 \geq [\ov{L}^2
- \log (D (g,m, m - g - 2) (m + g))[K:\Qb] ]
 / (2m \, [K:\Qb]) \, 
$$
(resp.
$$
\lb_{m-g} (\ov L) + 1 
 \geq [\ov{L}^2
- \log (D (g,m, m - g - 1) (m + g))[K:\Qb] ]
 / (2m \, [K:\Qb]) \,  . $$
}

\medskip

\noindent {\it Remark.} A lower bound $\lb_1 (\ov L)$ is given by 
Theorem~2 in [S]. However the statement of loc.cit. is not correct 
when $\deg_K (L) = m = 1$. In that case, we deduce from [S], 
Corollary 1, that
$$
{\ov L}^2 \leq 2 \, (\log \, \Vert e \Vert + 1) \, [K:\Qb] \, ,
$$
which is slightly weaker than [S], Th.~2.

\medskip

\noindent {\bf 3.7. Proof of Theorem 4.} Let $\pi : \wt X \ra X$ be 
a resolution of $X$. The canonical map
$$
H^1 (X , L^{-1}) \ra H^1 (\wt X , \pi^* (L^{-1}))
$$
is split injective, and it induces an isometry of hermitian complex 
vector spaces. Therefore, it is enough to prove Theorem 4  for 
$\wt X$ and $\pi^* (L^{-1})$, i.e. we can assume that $X$ is 
regular. (The same argument will prove Theorem 4 under the 
assumption that $X$ is a normal scheme, not necessarily semi-stable 
over $S$.)

\smallskip

If $\ov K$ is an algebraic closure of the field $K$,
the secant variety
$$
\Si_d \sbs \Pb (H^1 (X , L^{-1})_{\ov K})
$$
is defined over $K$ for any integer $d$. Its degree is $D(g,m,d)$.

\smallskip

For any $e \in H^1 (X,L^{-1})$ let
$$
\Vert e \Vert = \sup_{\s} \, \Vert e \Vert_{\s}
$$
where $\s$ runs over all complex embeddings $K \hra \Cb$, and $\Vert 
e \Vert_{\s}$ is the corresponding $L^2$-norm. Choose vectors $e_1 , 
\ldots , e_{m+g-1}$ in $H^1 (X , L^{-1})$ such that, for any $j = 1, 
\ldots , m + g - 1$, we have
$$
\Vert e_j \Vert = \lb_j (\ov L) \, .
$$
Given $k$ as in the statement of Theorem 4 i) or $k= m-g$
(resp. $k= m-g-1$) in case ii), we 
let $V$ be the $\Oc_K$-module spanned 
by $e_1 , e_2 , \ldots , e_k$. By Theorem 1 and Proposition 1,
 the projective space $\Pb (V_K)$ does not contain 
$\Si_{k-1}$. Using Proposition~4, we know that there exists an 
extension class $e \in H^1 (X , L^{-1})$ such that $[e] \notin 
\Si_{k-1}$ and
$$
\log \, \Vert e \Vert \leq \lb_{\max} (\ov V) + \log \, (D  
(g,m,k-1) (m+g)) \, , \leqno (15)
$$
with
$$
\lb_{\max} (\ov V) = \lb_k (\ov L) \, . \leqno (16)
$$

Since $e \notin \Si_{k-1}$ the corresponding extension $E$ of $L$ 
by $\Oc_X$ is either semi-stable when restricted to $X_{\ov K}$ or
, in case i), 
such that there exists a line bundle $M_K$ contained in $E_K$ with
$$
\deg(M_K) = m-d \, ,
$$
with $d \geq k$. When $E_{\ov K}$ is semi-stable, [S] Corollary 1 
asserts that
$$
{\ov L}^2 \leq (2m\log \, \Vert e \Vert + 2) [K:\Qb] \, . \leqno 
(17)
$$
In the remaining case, the inequality (35) in [S] reads
$$
d^2 \,  \ov{L}^2 + (m^2 - 2 \, m \, d) \, d \, e \, (\ov L) \leq m^2 \, 
[K:\Qb] \, d \, \log \, \Vert e \Vert + m^2 \, [K:\Qb] / 2 \, .
$$
Since $d \geq k $ and since
$$
2 \, m \, e \, (\ov L) \leq \ov{L}^2 \, \leqno (18)
$$
([Z] Theorem 6.3, valid under our assumption on $\ov L$)
we deduce that
$$
k (\ov{L}^2 - 2 \, m \, e \, (\ov L)) + m^2 e \, (\ov L) \leq m^2 \, 
[K:\Qb] \, (\log \, \Vert e \Vert + 1) \, . \leqno (19)
$$
From (15), (16) and (19) it follows that
$$
\lb_k (\ov L) + 1 \geq [k (\ov{L}^2 - 2 \, m \, e \, (\ov L)) + m^2 e 
\, (\ov L) - \log (D (g,m, k - 1)
(m + g))[K:\Qb] ] / (m^2 \, [K:\Qb]) \, 
. \leqno (20)
$$
On the other hand, (15), (16) and (17) imply
$$
\lb_k (\ov L) + 1 \geq [\ov{L}^2
- \log (D (g,m, k - 1) (m + g))[K:\Qb] ]
 / (2m \, [K:\Qb]) \, . \leqno (21)
$$
 This proves ii). In case i), we have
 to check that the right hand side
 in (21) is greater or equal to the right hand side
 in (20). But this follows easily
 from (18) and $m >2k$.
 This proves (15).

\smallskip

The inequality (16)  follows by Serre duality from (15) and 
Proposition~3. \hfill q.e.d.

\bigskip

\noindent {\bf 3.8.} Theorem~4 provides a precise lower (resp. 
upper) bound for the successive minima of $H^1 (X,L^{-1})$ (resp. 
$H^0 (X , L \ot \om_{X/S})$), since we can combine it with the lower 
bound for $e \, (\ov L)$ given in Prop.~4.

\smallskip

Let us make this more explicit when
$$
\ov L = \ov{\om}_{X/S}^{\ot n},
$$
where $n$ is a large positive integer,
hence $m = 2(g-1)n$. Assume $g > 1$. Proposition~4 
reads
$$
e \, (\ov L) \geq n \ {\ov{\om}_{X/S}^2 \over 4g \, (g-1)} \, .
$$
Therefore, since the coefficient of $ e \, (\ov L)$ is non-negative,
$$
k (\ov{L}^2 - 2 \, m \, e \, (\ov L)) + m^2 e \, (\ov L) \geq n^2 
(k + n) \, {g-1 \over g} \ \ov{\om}_{X/S}^2 \, .
$$
From Theorem~4 i) we get
$$
(\lb_k (\ov L) + 1) [K:\Qb] \geq {k + n \over 4 g (g-1)} \ 
\ov{\om}_{X/S}^2 - 
\log \, (D (g,2n (g-1),k-1) (m+g))[K:\Qb]/ m^2 
 \, . \leqno (22)
$$

\smallskip

We have
$$
D(g,2n (g-1), k-1) = \sum_{\a = 0}^g \left( \matrix{ m + g-1 
- k - 
\a \cr k - \a \cr} \right)
\left( \matrix{g \cr \a \cr} \right) 
\leq C_1(g) m! \, ,
$$
where $C_1(g)$ is a positive constant depending on
$g$ only. The Stirling formula implies
$$
\log \, D (g,2n (g-1) , k-1) \leq C_2(g)
m\log (m) .
$$
Therefore the second summand in the right hand side of (22)
is bounded above by a constant multiple
of $m \log (m)/ (2m^2) = \log(m)/2 $.

\smallskip

 Let $\om^2 = \ov{\om}_{X/S}^2 / [K:\Qb]$. We get the following:

\bigskip

\noindent {\bf Corollary 1.} {\it There exist constants
 $C(g) > 0$ such that,  if
$
\ov L = \ov{\om}_{X/S}^{\ot n}\,
 , n\geq 1$, $k < (g-1)n$
 and 
 $m = 2 (g-1) \, n$,
 the following inequality holds
 
 $$
\lb_k (\ov{\om}_{X/S}^{\ot n})  \geq {n + k  \over 4g (g-1)} 
\,   \om^2 \, - C(g) \log(n+1) \, .
$$
 
}

\bigskip

\noindent {\bf 3.9.}  From Theorem~4 one can also get an 
explicit lower bound for $\mu_k \, (\ov{\om}_{X/S}^{\ot 
n})$.

\smallskip

First, since $k \leq m/2$, we get, as in 3.8,
$$
\sum_{j=1}^k \log D (g,m,d(j)-1) \, {m+g \over m^2} \leq C_3(g)  \, n \log (n+1)  \, . 
\leqno (23)
$$
Similarly, by Stirling's formula and the known formula for the 
euclidean volume of the unit ball,
$$
C (m+g-2,K) = {m+g-1 \over 2} \, \log \vert \D_K \vert + O 
(m \log m) [K : \Qb] \, . \leqno (24)
$$

On the other hand, using Proposition~4 as in 3.8 above, we 
get
$$
{k (k+1) \over 2} (\ov{L}^2 - 2 \, m \, e 
(\ov L)) + k \, m^2 \, e (\ov L) \geq n^2 \ {g-1 \over g} 
\left[ {k (k+1) \over 2} + k \, n \right] \, 
\ov{\om}_{X/S}^2 \, . \leqno (25)
$$

If we put (14), (23), (24), (25) together we conclude that:

\bigskip

\noindent {\bf Corollary 2.} {\it There exist  constants $C(g) > 0$  such 
that, if $\ov L = \ov{\om}_{X/S}^{\ot n}$, $n \geq 1$, we 
have, when $k < (g-1)n$,
$$
\mu_k (\ov L) \geq  \, {1 \over 4g(g-1)} \left[ {k(k+1) \over 
2} + kn \right] \, \ov{\om}_{X/S}^2 - (2n+1) \, {g-1 
\over 2} \, \log \vert \D_K \vert - C(g) n \log (n+1) \, .
$$

}

\vglue 2cm

\noindent {\centerline {\bf References}}

\medskip

\item{[A]} A. Arakelov: Intersection theory of divisors on
an arithmetic surface, {\it Math. USSR, Izv. } {\bf 8},
1974, 1167-1180.

\item{[ACGH]} E. Arbarello, M. Cornalba, P.A. Griffiths, J. Harris : 
Geometry of Algebraic Curves, Vol.~I, 1985, Springer-Verlag.

\item{[B]} A. Bertram : Moduli of rank $2$ vector bundles, theta 
divisors, and the geometry of curves in projective space, {\it J. 
Diff. Geom.} {\bf 35}, 1992, 429-469.

\item{[BH]} R.C. Barker, G.Harman :
The difference between consecutive primes,
{\it Proc. Lond. Math. Soc.} {\bf 72}, 1996, 261-280.

\item{[Bo-Va]} E.Bombieri, J. Vaaler :
On Siegel's lemma, {\it Invent. Math.} {\bf 73}, 1983, 11-32.

\item{[BGS]} J.-B. Bost, H. Gillet, C. Soul\'e: Heights of projective 
varieties and positive Green forms, {\it J. Am. Math. Soc.}
{\bf 7}, 1994, 903-1027.

\item{[E]} R. Elkik : Fonctions de Green, Volumes de Faltings, 
Application aux surfaces arithm\'etiques, 
{\it Ast\'erisque}
{\bf 127}, 1985, 89-112.

 \item{[I]} A.Ivic : 
 The Riemann zeta-function. The theory of the Riemann
 zeta-function with applications, 
 A Wiley-Interscience Publication,
 New York etc., 1985, John Wiley \& Sons.

\item{[N]} J.Nagura : 
 On the interval containing at least one prime number,
 {\ it Proc. Japan Acad.} {\bf 28}, 1952, 177-181. 

 \item{[R]} P. Ribenboim :
 The book of prime number records, New York etc.: Springer-Verlag,
 1988 .
 
\item{[S]} C. Soul\'e : A vanishing theorem on arithmetic surfaces,
{\it Invent. Math.} {\bf 116}, 1994, 577-599.

\item{[V]} C. Voisin : Appendix: On linear subspaces contained in the 
secant varieties of a projective curve. This volume.

\item{[Z]} S. Zhang : Positive line bundles on arithmetic surfaces,
{\it Annals of Maths.} {\bf 136}, 1992, 569-587.
\bye